\documentclass[12pt]{article}
\usepackage[margin=2cm]{geometry}
\RequirePackage{setspace}
\usepackage{amssymb,amsthm,amsmath, tikz,mathtools,subfig, mathrsfs,array}
\usetikzlibrary{shapes}
\newtheorem{thm}{Theorem}[section]

\newtheorem{lemma}[thm]{Lemma}

\newtheorem{cor}[thm]{Corollary}

\begin{document}

\title{Enumerating the Digitally Convex Sets of Powers of Cycles and Cartesian Products of Paths and Complete Graphs}
\author{MacKenzie Carr \footnote{Department of Mathematics and Statistics, University of Victoria, Victoria, BC, Canada}  \\ \small{\texttt{mackenziecarr@uvic.ca}} \and Christina M.~Mynhardt \footnotemark[1] \thanks{Supported by an NSERC Grant CANADA, Grant number  RGPIN-5442-2015} \\ \small{\texttt{kieka@uvic.ca}} \and Ortrud R.~Oellermann \footnote{Department of Mathematics and Statistics, University of Winnipeg, Winnipeg, MB, Canada} \thanks{Supported by an NSERC Grant CANADA, Grant number RGPIN-2016-05237} \\ \small{\texttt{o.oellermann@uwinnipeg.ca}}}

\date{\vspace{-5ex}}

\maketitle

\begin{abstract}
Given a finite set $V$, a \emph{convexity} $\mathscr{C}$, is a collection of subsets of $V$ that contains both the empty set and the set $V$ and is closed under intersections. The elements of $\mathscr{C}$ are called \emph{convex sets}. The digital convexity, originally proposed as a tool for processing digital images, is defined as follows: a subset $S\subseteq V(G)$ is \emph{digitally convex} if, for every $v\in V(G)$, we have $N[v]\subseteq N[S]$ implies $v\in S$. The number of cyclic binary strings with blocks of length at least $k$ is expressed as a linear recurrence relation for $k\geq 2$. A bijection is established between these cyclic binary strings and the digitally convex sets of the $(k-1)^{th}$ power of a cycle. A closed formula for the number of digitally convex sets of the Cartesian product of two complete graphs is derived. A bijection is established between the digitally convex sets of the Cartesian product of two paths, $P_n \square P_m$, and certain types of $n \times m$ binary arrays. 
\end{abstract}

\section{Introduction}

Given a finite set $V$, a collection $\mathscr{C}$, of subsets of $V$ is called a \emph{convexity} or \emph{alignment} if it contains $\emptyset$ and $V$ and is closed under intersections. The elements of a convexity $\mathscr{C}$ are called \emph{convex sets} and the ordered pair $(V,\mathscr{C})$ is an \emph{aligned space}. For any subset $S\subseteq V$, the \emph{convex hull} of $S$, denoted by $CH_\mathscr{C}(S)$, is the smallest convex set that contains $S$. For any $S\subseteq V$, if $CH_\mathscr{C}(S)=S$, then $S$ is a convex set. Van de Vel provides an in-depth study of abstract convex structures in~\cite{theory}. 

There are several convexities defined on the vertex set of a graph. The most natural extension of the Euclidean convexity to graphs is defined using an interval notion. For $a,b\in V(G)$, the collection of vertices that are on some $a$-$b$ geodesic (shortest $a$-$b$ path) forms the \emph{geodesic interval} between $a$ and $b$. Then, a set $S\subseteq V(G)$ is \emph{g-convex} if it contains the geodesic interval between every pair of vertices in $S$. The collection of all $g$-convex sets in a graph $G$ is a convexity called the \emph{geodesic convexity} of $G$. 

Several other graph convexities defined in terms of different types of intervals between pairs of vertices were studied, for example, in \cite{triangle,hhd,hyper}. Interval structures between three or more vertices have led to yet more graph convexities examined, for example, in \cite{3steinsimp, mintree, steiner}. 

In this paper we study the digital convexity of a graph, introduced by Rosenfeld and Pfaltz in~\cite{seq} as a tool for processing digital images. Rather than using a definition based on an interval structure, the digital convexity is instead defined in terms of neighbourhoods. The \emph{open neighbourhood} of a vertex $v\in V(G)$, denoted by $N_G(v)$ or $N(v)$ when the graph $G$ is obvious, is defined as $N_G(v)=\{x\in V(G)\mid xv\in E(G)\}$. Similarly, the \emph{closed neighbourhood} of $v$, denoted by $N_G[v]$ or $N[v]$, is defined as $N_G[v]=N_G(v)\cup\{v\}$. For a set $S\subseteq V(G)$, the closed neighbourhood of $S$, denoted by $N_G[S]$ or $N[S]$, is defined as $N_G[S]=\bigcup_{v\in S}N_G[v]$. 

\begin{figure}[ht]
\centering
\subfloat[]{
\begin{tikzpicture}[thick]

\foreach \x/\y in {1/1, 0.5/1, 1/2, 1.5/2.5, 2/1.5, 2/2, 2.5/1.5} {
    \node [fill=black, draw=none, thick, minimum size=0.5cm] 
      at (\x-.25,2.5+0.25-\y) {};
    }
    \draw[step=0.5] (0, 0) grid (2.5,2.5);

\end{tikzpicture}
\label{fig:imageorig}
}
\hspace{15mm}
\subfloat[]{
\begin{tikzpicture}[thick]

\foreach \x/\y in {0.5/0.5, 0.5/1.5, 1/1.5, 1.5/1.5, 1.5/2, 1/1, 0.5/1, 1/2, 1.5/2.5, 2/1.5, 2/2, 2.5/1.5} {
    \node [fill=black, draw=none, thick, minimum size=0.5cm] 
      at (\x-.25,2.5+0.25-\y) {};
    }
    \draw[step=0.5] (0, 0) grid (2.5,2.5);
   
\end{tikzpicture}
\label{fig:imagesmooth}
}
\caption{A black and white digital image in (a) and its corresponding digital convex hull in (b)}
\label{fig:image}
\end{figure}
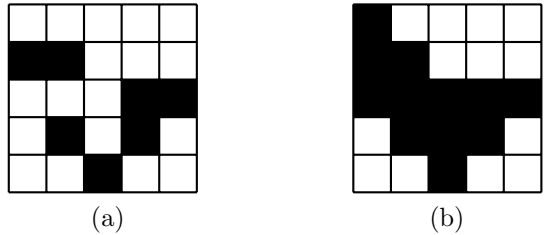

A set $S\subseteq V(G)$ is \emph{digitally convex} if $N_G[v]\subseteq N_G[S]$ implies $v\in S$ for every $v\in V(G)$. For a vertex $v\in V(G)$ and a set $S\subseteq V(G)$, if $N_G[v]-N_G[S-\{v\}]\neq \emptyset$, we say that $v$ has a \emph{private neighbour with respect to} $S$ in $G$. Thus, $S$ is digitally convex if and only if, for every $v\not\in S$, $v$ has a private neighbour with respect to $S$. Note that private neighbours are not necessarily unique and a vertex $v$ can be a private neighbour for multiple vertices. For a graph $G$, the collection of all digitally convex sets in $G$ is the \emph{digital convexity} of $G$, denoted by $\mathscr{D}(G)$. The number of digitally convex sets in $G$ is denoted by $n_\mathscr{D}(G)$. The digital convexity is studied in the context of closure systems in~\cite{closure}. The relationship between digital convexity and domination in a graph is examined in~\cite{frame,digconv}. In particular, for $v\in V(G)$ and $S\subseteq V(G)$, if $N[v]\subseteq N[S]$, then $S$ is a \emph{local dominating set for} $v$. Thus, a digitally convex set is a set of vertices containing every vertex for which it is a local dominating set. Figure~\ref{fig:image}(a) shows a black and white digital image (or grid), with the digital convex hull of the black pixels shown in Figure~\ref{fig:image}(b). 

The problem of determining the number of convex sets of a graph has been studied for several of the above graph convexities. In the case of the geodesic convexity, it has been shown that the number of $g$-convex sets of a tree is equal to the number of its subtrees, a problem which is explored in~\cite{subtrees,binary,enumsub}. Brown and Oellermann~\cite{spectrum} determined that the problem of enumerating the $g$-convex sets of a cograph can be performed in linear time but, for an arbitrary graph, the problem is \#P-complete. Graphs with a minimal number of $g$-convex sets or a minimal number of $m$-convex sets was examined in~\cite{minimalconv}. An algorithm for generating the digitally convex sets of a tree, as well as sharp upper and lower bounds on the number of digitally convex sets of a tree and a closed formula for the number of digitally convex sets of a path are given in~\cite{enumtree}.   

This paper focuses on the enumeration of digitally convex sets in powers of cycles and in Cartesian products of two complete graphs and of two paths. In Section 2, we derive a recurrence relation for cyclic binary strings, all of whose blocks have length at least $k \geq 2$. We then establish a bijection between the digitally convex sets of the $k^{th}$ power of a cycle and the cyclic binary strings whose blocks all have length at least $k+1$. In Section 3, we develop a closed formula for the number of digitally convex sets of the Cartesian product of two complete graphs, $K_n\square K_m$. It is shown that there is a bijection between the number of digitally convex sets of the Cartesian product of two paths, $P_n\square P_m$, and certain types of $n\times m$ arrays. In the special case where $m=2$, the number of digitally convex sets in the Cartesian product $P_n\square P_2$ can be expressed as a linear recurrence relation of order three.

\section{Digital Convexity in Powers of Cycles}

In this section, we establish a recurrence relation for the number of cyclic binary strings with blocks of length at least $k\geq 2$ and we show that the number of digitally convex sets in the $(k-1)^{th}$ power of a cycle satisfies this same recurrence. Throughout this section, for a cycle $C_n$, we denote the vertices by $v_1,v_2,\dots,v_n$ with $v_iv_{i+1}\in E(C_n)$ for $i=1,2,\dots,n-1$ and $v_1v_n\in E(C_n)$. Recall that for a positive integer $d$, the $d^{th}$ power of a graph $G$ is the graph $G^d=(V,E')$, such that $uv\in E'$ if and only if $u$ and $v$ are distance at most $d$ apart in $G$. Note that for $3\leq n \leq 2k+1$, the graph $C^k_n$ is a complete graph. 

For $k\geq 2$, we let $\mathscr{B}_{k,n}$ be the set of cyclic binary strings of length $n$ in which each block (maximal run of 0's or 1's) has length at least $k$ if $n\geq k$, or length exactly $n$ if $n<k$. Let $a_k(n)=|\mathscr{B}_{k,n}|$. Now, we establish a recurrence relation satisfied by $a_k(n)$. Munarini and Salvi use the Sch\"utzenberger symbolic method to show this relation for $k=2$ in~\cite{zigzag}. We use this same method to generalize the result of Munarini and Salvi to any $k\geq 2$. 

\begin{lemma}
Let $k\geq 2$. Then $a_k(i)=2$ for $3\leq i\leq 2k-1$, $a_k(j)=2+j(j-2k+1)$ for $2k\leq j\leq 2k+2$ and, for $n\geq 2k+3$, \begin{center}$a_k(n)=2a_k(n-1)-a_k(n-2)+a_k(n-2k).$\end{center} 
\label{lem:bin}
\end{lemma}

\noindent
\emph{Proof.} First, we establish the initial conditions. If $3\leq n < k$, then the only cyclic binary strings of length $n$ with each block of length exactly $n$ are clearly $(00\dots0)$ and $(11\dots1)$. Thus, $a_k(n)=2$ for $3\leq n<k$. Now suppose $k\leq n\leq 2k-1$. Any cyclic binary string with at least two blocks must have one block of length $\ell\leq n/2 \leq (2k-1)/2<k$. So the only strings in $\mathscr{B}_{k,n}$ are $(00\dots0)$ and $(11\dots1)$ when $k\leq n\leq 2k-1$. 

If $2k\leq n\leq 2k+2$, then clearly both $(0 0\dots0)$ and $(1 1\dots1)$ are cyclic binary strings in $\mathscr{B}_{k,n}$. The remaining strings in $\mathscr{B}_{k,n}$ are those with two blocks, one of length $\ell$, with $k\leq \ell \leq n-k$, and the other of length $n-\ell$. Without loss of generality, let the block of length $\ell$ be $\ell$ consecutive 1's. There are $n$ distinct cyclic shifts of these two blocks, giving $n$ distinct cyclic binary strings with $\ell$ consecutive 1's. There are $n-2k+1$ possible values of $\ell$, so there are $n(n-2k+1)$ cyclic binary strings in $\mathscr{B}_{k,n}$ with exactly two blocks. Overall, we have $a_k(n)=2+n(n-2k+1)$ for $2k\leq n\leq 2k+2$. Therefore, the initial conditions hold. 

Now, we find the generating function to show the desired recurrence.  To do this, we uniquely decompose the strings in $\mathscr{B}_{k,n}$ into the smaller strings $(0\dots0)$, $(1\dots1)$, $(01)$, and $(10)$. The latter two types of strings will be called \emph{principal blocks}. Then, the strings in $\mathscr{B}_{k,n}$ containing a principal block can be decomposed in one of the following ways:
$$(0\dots0)(01)(1^{k_1})(10)(0^{k_2})\dots$$
$$(1\dots1)(10)(0^{k_1})(01)(1^{k_2})\dots$$
where $k_i\geq k-2$ and an exponent of $k_i$ means that the indicated bit is repeated a total of $k_i$ times. So each principal block $(01)$ that does not appear at the beginning of the string must be preceded by a string of type $(0^{k_i})$ and followed by a string of type $(1^{k_{i+1}})$, and, similarly, each principal block $(10)$ must be preceded by a string of type $(1^{k_j})$ and followed by a string of type $(0^{k_{j+1}})$. 
Now, we can break these two cases down further into those containing an even number of principal blocks and those containing an odd number. The strings in $\mathscr{B}_{k,n}$ containing $2\ell>0$ principal blocks have one of the forms
$$(0^{p_1})(01)(1^{k_1})(10)(0^{k_2})\dots(01)(1^{k_{2\ell-1}})(10)(0^{p_2})$$
$$(1^{q_1})(10)(0^{k_1})(01)(1^{k_2})\dots(10)(0^{k_{2\ell-1}})(01)(1^{q_2})$$
where each $k_i \geq k-2$, $p_1+p_2 \geq k-2$ and $q_1+q_2 \geq k-2$. 

Similarly, the strings in $\mathscr{B}_{k,n}$ containing $2\ell+1$ principal blocks have one of the forms 
$$(0^{r_1})(01)(1^{k_1})(10)(0^{k_2})\dots(01)(1^{r_2})$$
$$(1^{s_1})(10)(0^{k_1})(01)(1^{k_2})\dots(10)(0^{s_2})$$ 
where, as above, each $k_i \geq k-2$, and $r_1,r_2,s_1,s_2 \geq k-1$.

Then, we can express $\mathscr{B}_k$, the set of all cyclic binary strings with all blocks of length at least $k$, using the symbolic method. First, we let $0^*=\{\varepsilon, 0, 00, \dots\}$ ($\varepsilon$ denotes the empty string), $1^*=\{\varepsilon, 1, 11, \dots\}$, $0^+=0^* - \{\varepsilon\}$ and $1^+=1^* - \{\varepsilon\}$. Then, 

\begin{equation*}\begin{gathered}\mathscr{B}_k = 0^+\cup1^+\bigcup_{p=0}^{k-2}\big(\bigcup_{\ell=1}^{\infty} \underbrace{0^p(01)1^{k-2}1^*(10)0^{k-2}0^*\dots(01)1^{k-2}1^*(10)0^{k-2-p}0^*}_{2\ell~principal~blocks} \big)\\
\bigcup_{\ell=1}^{\infty} \underbrace{0^*0^{k-2}(01)1^{k-2}1^*(10)0^{k-2}0^*\dots(01)1^{k-2}1^*(10)0^*}_{2\ell~principal~blocks} \\
\bigcup_{q=0}^{k-2}\big(\bigcup_{\ell=1}^{\infty} \underbrace{1^q(10)0^{k-2}0^*(01)1^{k-2}1^*\dots(10)0^{k-2}0^*(01)1^{k-2-q}1^*}_{2\ell~principal~blocks} \big)\\
\bigcup_{\ell=1}^{\infty} \underbrace{1^*1^{k-2}(10)0^{k-2}0^*(01)1^{k-2}1^*\dots(10)0^{k-2}0^*(01)1^*}_{2\ell~principal~blocks}\\
\bigcup_{\ell=0}^{\infty}\underbrace{0^*0^{k-1}(01)1^{k-2}1^*(10)0^{k-2}0^*\dots(01)1^{k-1}1^*}_{2\ell+1~principal~blocks} \\
\bigcup_{\ell=0}^{\infty}\underbrace{1^*1^{k-1}(10)0^{k-2}0^*(01)1^{k-2}1^*\dots(10)0^{k-1}0^*}_{2\ell+1~principal~blocks}\end{gathered}\end{equation*} 
Note that we divide the strings with $2\ell$ principal blocks into two cases: those beginning with fewer than $k-1$ 0's (or $k-1$ 1's) and those beginning with at least $k-1$ 0's (resp.~1's). We divide in this way because, in the first case, there is a minimum number of 0's (resp.~1's) that must be at the end of the string so that it is contained in $\mathscr{B}_{k,n}$. There is no such minimum in the second case. 

Now, we can use the symbolic method to find the generating function for $a_k(n)$. The generating function for a string $0^*$ or $1^*$ is $1/(1-x)$, the generating function for a string $0^i$ or $1^{i}$ is $x^{i}$, and the generating function for a principal block $(01)$ or $(10)$ is $x^{2}$. We multiply the generating functions for the substrings to get a generating function for $a_k(n)$, using the same deconstruction as above. 
\begin{equation*}\begin{gathered}B(x) = \sum_{n=0}^{\infty} a_k(n)x^n = 2\frac{x}{1-x} 
\end{gathered}\end{equation*}
\begin{equation*}\begin{gathered} + \sum_{p=0}^{k-2} \sum_{\ell=1}^{\infty}\bigg(x^p \frac{x^{k-2-p}}{1-x} (x^2)^{2\ell} \frac{(x^{k-2})^{2\ell-1}}{(1-x)^{2\ell-1}}\bigg)  +\sum_{\ell=1}^{\infty}\bigg(\frac{x^{k-1}}{1-x} \frac{1}{1-x} (x^2)^{2\ell} \frac{(x^{k-2})^{2\ell-1}}{(1-x)^{2\ell-1}}\bigg)
\end{gathered}\end{equation*}
\begin{equation*}\begin{gathered}
 + \sum_{q=0}^{k-2} \sum_{\ell=1}^{\infty}\bigg(x^q \frac{x^{k-2-q}}{1-x} (x^2)^{2\ell} \frac{(x^{k-2})^{2\ell-1}}{(1-x)^{2\ell-1}}\bigg)  +\sum_{\ell=1}^{\infty}\bigg(\frac{x^{k-1}}{1-x} \frac{1}{1-x} (x^2)^{2\ell} \frac{(x^{k-2})^{2\ell-1}}{(1-x)^{2\ell-1}}\bigg) 
 \end{gathered}\end{equation*}
\begin{equation*}\begin{gathered} + \sum_{\ell=0}^{\infty} \bigg( (x^2)^{2\ell+1} \frac{(x^{k-2})^{2\ell}}{(1-x)^{2\ell}} \frac{(x^{k-1})^2}{(1-x)^2}\bigg) + \sum_{\ell=0}^{\infty} \bigg((x^2)^{2\ell+1} \frac{(x^{k-2})^{2\ell}}{(1-x)^{2\ell}} \frac{(x^{k-1})^2}{(1-x)^2}\bigg)\end{gathered}\end{equation*}
\begin{equation*}\begin{gathered}
 = \frac{2x}{1-x} + 2\sum_{p=0}^{k-2} \sum_{\ell=1}^{\infty} \bigg(\frac{x^{2k}}{(1-x)^2}\bigg)^{\ell} + \frac{2x}{1-x} \sum_{\ell=1}^{\infty} \bigg(\frac{x^{2k}}{(1-x)^2}\bigg)^{\ell} + \frac{2x^{2k}}{(1-x)^2} \sum_{\ell=0}^{\infty} \bigg(\frac{x^{2k}}{(1-x)^2}\bigg)^{\ell}
\end{gathered}\end{equation*}

Each series in terms of $\ell$ in $B(x)$ is a geometric series. So we can write the function in closed form and simplify.

\begin{equation*}\begin{gathered}B(x) = \frac{2x}{1-x} + 2\sum_{p=0}^{k-2} \bigg(\frac{\frac{x^{2k}}{(1-x)^2}}{1-\frac{x^{2k}}{(1-x)^2}}\bigg) + \frac{2x}{1-x}\bigg(\frac{\frac{x^{2k}}{(1-x)^2}}{1-\frac{x^{2k}}{(1-x)^2}}\bigg) + \frac{2x^{2k}}{(1-x)^2} \bigg(\frac{1}{1-\frac{x^2k}{(1-x)^2}}\bigg)\end{gathered}\end{equation*}
\begin{equation*}\begin{gathered} = \frac{2x}{1-x} + \sum_{p=0}^{k-2}\bigg(\frac{2x^{2k}}{1-2x+x^2-x^{2k}}\bigg) + \frac{2x^{2k+1}}{(1-x)(1-2x+x^2-x^{2k})} + \frac{2x^{2k}}{1-2x+x^2-x^{2k}}\end{gathered}\end{equation*}
\begin{equation*}\begin{gathered} = \frac{2x}{1-x} + \frac{2(k-1)x^{2k}}{1-2x+x^2-x^{2k}} + \frac{2x^{2k+1}}{(1-x)(1-2x+x^2-x^{2k})} + \frac{2x^{2k}}{1-2x+x^2-x^{2k}}\end{gathered}\end{equation*}
\begin{equation*}\begin{gathered} = \frac{2x-4x^2+2x^3+2kx^{2k}-2kx^{2k+1}}{(1-x)(1-2x+x^2-x^{2k})}\end{gathered}\end{equation*}
\begin{equation*}\begin{gathered} = \frac{2x-2x^2+2kx^{2k}}{1-2x+x^2-x^{2k}}.
\end{gathered}\end{equation*}

From the form of the generating function, we know that $a_k(n)-2a_k(n-1)+a_k(n-2)-a_k(n-2k)=0.$ Rearranging this, we get the desired recurrence. 
\qed 

By establishing a suitable bijection, we now show that the recurrence relation established in Lemma~\ref{lem:bin} holds for the number of digitally convex sets of powers of cycles. 

\begin{thm}
Let $C^k_n$ be the $k^{th}$ power of the cycle $C_n$, with $k\geq 1$. Then $n_\mathscr{D}(C^k_i) = 2$ for $3\leq i \leq 2k+1$, $n_\mathscr{D}(C^k_j) = 2+j(j-2k-1)$ for $2k+2 \leq j \leq 2k+4$ and, for $n\geq 2k+5$, 
\begin{center}$n_\mathscr{D}(C^k_n) = 2n_\mathscr{D}(C^k_{n-1})-n_\mathscr{D}(C^k_{n-2}) + n_\mathscr{D}(C^k_{n-2k-2}).$\end{center}
\label{thm:cyclepower}
\end{thm}

\noindent
\emph{Proof.} To prove the recurrence, we show a bijection between the digitally convex sets in $\mathscr{D}(C^k_n)$ and the cyclic binary $n$-bit strings in $\mathscr{B}_{k+1,n}$. If $n< k+1$, then these are the cyclic binary strings with blocks of length exactly $n$, i.e.~$(0\dots0)$ and $(1\dots1)$. Clearly, since $k+1\leq 2k+1$, there are exactly two digitally convex sets in $C^k_{n}$ when $n < k+1$, the sets $\emptyset$ and $V(C^k_n)$. These sets get mapped to $(0\dots0)$ and $(1\dots1)$, respectively. 

\begin{figure}[ht]
\centering 
\begin{tikzpicture}[thick, every node/.style={circle, draw=black, fill=black, inner sep=3}]

\node[label=left:$v_5$] (n1) at (0.25,0){};
\node[label=right:$v_4$] (n2) at (2.25,0){};
\node[label=right:$v_3$] (n3) at (3,1.25){};
\node[label=right:$v_2$] (n4) at (2.75,2.5){};
\node[inner sep = 3.5,fill=none,label=$v_1$] (n5) at (1.25,3.25){};
\node[inner sep = 3.5,fill=none, label=left:$v_7$] (n6) at (-0.25,2.5){};
\node[label=left:$v_6$] (n7) at (-0.5,1.25){};

\draw (n1) -- (n2);
\draw (n2) -- (n3);
\draw (n3) -- (n4);
\draw (n4) -- (n5);
\draw (n5) -- (n6);
\draw (n6) -- (n7);
\draw (n7) -- (n1);
\draw (n1) -- (n3);
\draw (n2) -- (n4);
\draw (n3) -- (n5);
\draw (n4) -- (n6);
\draw (n5) -- (n7);
\draw (n6) -- (n1);
\draw (n7) -- (n2);

\end{tikzpicture}
\caption{The digitally convex set $S=\{v_1, v_7\}$ of $C^2_7$ is indicated by white vertices}
\label{fig:c27}
\end{figure}
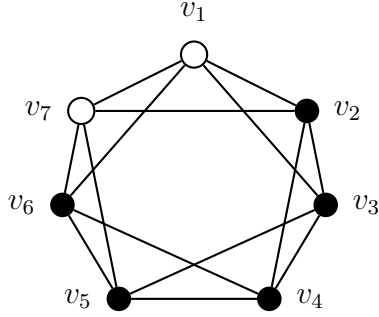

Now, suppose $n\geq k+1$. Then $\mathscr{B}_{k+1,n}$ is the set of cyclic binary $n$-bit strings whose maximal blocks each have length at least $k+1$. Given a digitally convex set $S\in\mathscr{D}(C^k_n)$, we get a corresponding cyclic binary $n$-bit string $S^*$ in the following way. For each vertex $v_i\in S$, set bits $i, i+1, \dots, i+k$ in $S^*$ to be 1, taking the index mod $n$ if $i+j>n$. After repeating this for each vertex in $S$, set the remaining bits in $S^*$ to 0. As an example, the digitally convex set $S=\{v_1,v_7\}$ in $C^2_7$, shown in Figure~\ref{fig:c27}, corresponds to the cyclic binary string $S^*=(1110001)$. 

It is clear from the construction of $S^*$ that each block of 1's must have length at least $k+1$. We show now that each block of 0's in $S^*$ must also have length at least $k+1$. Suppose that there is a block of 0's with length $\ell\leq k$, say bits $i, i+1, \dots, i+\ell-1$. Then, the vertex $v_i\not\in S$, since bit $i$ is 0 in $S^*$ and $v_{i+\ell}\in S$, since bit $i+\ell$ is 1 in $S^*$. In $C^k_n$, we must have $v_i v_{i+\ell}\in E(C^k_n)$, as $\ell \leq k$, as well as $v_{i+j}v_{i+\ell}\in E(C^k_n)$ for $j=1,2,\dots,\ell-1$.   In $S^*$, bit $i-1$ is also 1. So, by the construction of $S^*$, the vertex $v_{i-k-1}\in S$. This vertex is adjacent to the vertices $v_{i-k}, v_{i-k+1}, \dots,v_{i-1}$ in $C^k_n$. Thus, $N[v_{i}] \subseteq N[\{v_{i+\ell}, v_{i-k-1}\}] \subseteq N[S]$, contradicting the fact that $S$ is digitally convex in $C^k_n$. 

We now show that there is an injective map from $\mathscr{B}_{k+1,n}$ to the set of digitally convex sets of $C^k_n$. Let $S^*\in\mathscr{B}_{k+1,n}$. 
If $S^* = (00\dots0)$, then let $S=\emptyset$. If $S^* = (11\dots1)$, then let $S=V(C^k_n)$. Both of these are clearly digitally convex. Otherwise, let $B_1,B_2,\dots,B_r$ be the distinct blocks of at least $k+1$ 1's in $S^*$. Suppose bits $i, i+1, \dots, i+k+\ell-1$ are the bits of $B_1$. Then, let $S_1 = \{v_i, v_{i+1}, \dots, v_{i+\ell-1}\}$. Define $S_2,S_3,\dots,S_r$ similarly. Finally, let $S=S_1\cup S_2\cup\dots\cup S_r$. It is clear that $S$ would be mapped to $S^*$ using the above mapping. For example, if $S^* = (1110001)$ and $k=2$, then bits $7, 1, 2, 3$ are the bits of the only block of 1's. This string would be mapped to the set of vertices $\{v_7, v_1\}$, reversing the example of the mapping shown earlier in the proof. 
We show now that each such set $S$ is a digitally convex set in $\mathscr{D}(C^k_n)$. Suppose otherwise, i.e.~that $S$ is not digitally convex. Then, there must be some $v_j\not\in S$ such that $N[v_j] = \{v_{j-k}, v_{j-k+1}, \dots, v_j, v_{j+1}, \dots, v_{j+k}\} \subseteq N[S]$. 

Since $v_j\not\in S$ and $v_{j-k}\in N[S]$, we must have one of the vertices in $N[v_{j-k}]-\{v_j\}$ in $S$, say $v_{j-k+m}$ for some $m\in\{-k,-k+1,\dots,k-1\}$. Then, by definition of $S$, the bits $j-k+m$, $j-k+m+1,\dots, j+m$ are all 1 in $S^*$. None of these possible vertices $v_{j-k+m}$ is adjacent to $v_{j+k}$, so one of the vertices in $N[v_{j+k}]-\{v_j\}$ is in $S$, say $v_{j+k+p}$, for some $p\in\{-k+1, -k+2,\dots,k\}$. So, again by definition, the bits $j+k+p$, $j+k+p+1, \dots, j+2k+p$ are each 1 in $S^*$. In addition, these two vertices can be chosen so that each of the vertices $v_{j-k}, v_{j-k+1}, \dots, v_j, v_{j+1}, \dots, v_{j+k}$ appears in the closed neighbourhood of $v_{j-k+m}$ or $v_{j+k+p}$. Then, the maximum possible difference between $j-k+m$ and $j+k+p$ is $2k+1$. So the longest block of 0's in $S^*$ between bits $j+m$ and $j+k+p$ has length at most $k$, contradicting the fact that $S^*\in \mathscr{B}_{k+1,n}$. Therefore, $S$ is digitally convex in $C^k_n$. 

We have now shown a bijection between the digitally convex sets in $\mathscr{D}(C^k_n)$ and the cyclic binary strings in $\mathscr{B}_{k+1,n}$. So they satisfy the same recurrence. Therefore, $n_\mathscr{D}(C^k_n) = 2n_\mathscr{D}(C^k_{n-1})-n_\mathscr{D}(C^k_{n-2}) + n_\mathscr{D}(C^k_{n-2k-2)})$, with $n_\mathscr{D}(C^k_i) = 2$, for $3\leq i \leq 2k+1$, and $n_\mathscr{D}(C^k_j) = 2+j(j-2k-1)$, for $2k+2 \leq j \leq 2k+4$. \qed

\begin{cor}
Let $C_n$ be the cycle of order $n$. Then $n_\mathscr{D}(C_3)=2,~n_\mathscr{D}(C_4)=6,~n_\mathscr{D}(C_5)=12$, $ n_\mathscr{D}(C_6)=20$ and, for $n\geq 7$, \begin{center}$n_\mathscr{D}(C_n)=2n_\mathscr{D}(C_{n-1})-n_\mathscr{D}(C_{n-2})+n_\mathscr{D}(C_{n-4}).$\end{center}
\end{cor}

As shown above, the recurrence satisfied by $n_\mathscr{D}(C_n)$ is the same recurrence satisfied by the cyclic binary $n$-bit strings whose blocks each have length at least 2. This recurrence was shown, above and in~\cite{zigzag}, to have the generating function $$\frac{2x-2x^2+4x^4}{1-2x+x^2-x^4}.$$

Notice that this expands to 
\begin{equation*}\begin{gathered}2x + 2 x^2 + 2 x^3 + 6 x^4 + 12 x^5 + 20 x^6 + \dots + 74 x^9 +  122 x^{10} + 200 x^{11} \\ + \dots + 842 x^{14} + 1362 x^{15} + \dots + 9350 x^{19} + 15 126 x^{20} + \dots + 64080 x^{23} + 103 684 x^{24} \\+ 167 762 x^{25} + \dots + 710646 x^{28} + 1 149 852 x^{29} + 1 860 500 x^{30} + \dots\end{gathered}\end{equation*}
indicating that $C_5$ is the smallest cycle with more than ten digitally convex sets, $C_{10}$ is the smallest cycle with more than 100 digitally convex sets, $C_{15}$ is the smallest cycle with more than 1000 digitally convex sets, and $C_{20}$ is the smallest cycle with more than 10 000 digitally convex sets. However, $C_{24}$, not $C_{25}$, is the smallest cycle with more than 100 000 digitally convex sets and, similarly, $C_{29}$ is the smallest cycle with more than 1 000 000 digitally convex sets. This pattern suggests that $n_\mathscr{D}(C_{5k}) \geq 10^k$, but not that $C_{5k}$ is the smallest such cycle satisfying this inequality.

\section{Digital Convexity and Cartesian Products}

\begin{figure}[ht]
\centering
\begin{tikzpicture}[thick, every node/.style={circle, draw=black, fill=black, inner sep=3}]

\node[label=left:{$(1,1)$}] (n11) at (0,0){};
\node[fill=none] (n21) at (1.5,0){};
\node[label=right:{$(3,1)$}] (n31) at (3,0){};
\node[label=left:{$(1,2)$}] (n12) at (0,1.5){};
\node (n22) at (1.5,1.5){};
\node[label=right:{$(3,2)$}] (n32) at (3,1.5){};
\node[draw=none,fill=none,label=below:{$(2,1)$}] (n4) at (1.5,0.35){};
\node[draw=none,fill=none,label=above:{$(2,2)$}] (n4) at (1.5,1.15){};

\draw (n11) -- (n21);
\draw (n21) -- (n31);
\draw (n12) -- (n22);
\draw (n22) -- (n32);
\draw (n11) -- (n12);
\draw (n21) -- (n22);
\draw (n31) -- (n32);
\draw (n11) [looseness = 1.75, out = -35, in = -145] to (n31);
\draw (n12) [looseness = 1.75, out = 35, in = 145] to (n32);

\end{tikzpicture}

\caption{The set $\{(2,1)\}\in \mathscr{D}(K_3 \square K_2)$ is indicated in white}
\label{fig:k3k2}
\end{figure}
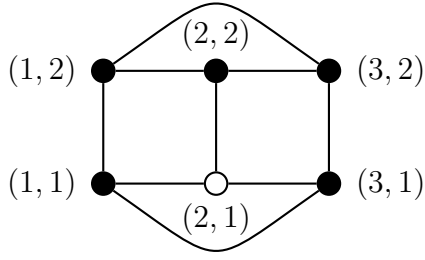

A digitally convex set in the Cartesian product $G\square H$ is not necessarily digitally convex when restricted to $G$ or to $H$. In other words, if $S\in \mathscr{D}(G\square H)$, then the set $S_G = \{x\in V(G) \mid (x,y)\in S\}$ is not necessarily digitally convex in $G$. As an example, the set $\{(2,1)\}$, shown in Figure~\ref{fig:k3k2}, is digitally convex in $K_3 \square K_2$ but $\{2\}\not\in\mathscr{D}(K_3)$ and $\{1\}\not\in\mathscr{D}(K_2)$, as the only digitally convex sets in a complete graph are the empty set and the entire vertex set. This example shows that, even in small graphs, the number of digitally convex sets of a Cartesian product of graphs $G$ and $H$ cannot be computed from those of $G$ and $H$ in an obvious manner. We begin by examining the number of digitally convex sets in the Cartesian product of complete graphs, $K_n \square K_m$, to show how different this number is from the number of digitally convex sets in either of the constituent graphs of the product. 

\begin{thm}
For any $m,n\geq 1$, $n_\mathscr{D}(K_n \square K_m) = 2+(2^n-2)(2^m-2)$.
\end{thm}

\noindent
\emph{Proof.} We begin by denoting the vertices of $K_n \square K_m$ by $(v_i,u_j)$ for $i=1,2,\dots,n$ and $j=1,2,\dots,m$. 

Now, let $\emptyset\neq S_1 \subsetneqq \{v_1,v_2,\dots,v_n\}$ and $\emptyset \neq S_2 \subsetneqq \{u_1,u_2,\dots,u_m\}$. Then, $S=S_1\times S_2 = \{(v_i,u_j)\mid v_i\in S_1,u_j\in S_2\}$ is digitally convex in $K_n\square K_m$. Consider $(v_x,u_y)\not\in S$. If $v_x\not\in S_1$ and $u_y\not\in S_2$, then $(v_x,u_y)\not\in N[S]$. If $v_x\in S_1$ and $u_y\not\in S_2$, then there is some $v_z\not\in S_1$ such that $(v_x,u_y)(v_z,u_y)\in E(K_n\square K_m)$ and $(v_z,u_y)\not\in N[S]$. If $v_x\not\in S_1$ and $u_y\in S_2$, then there is some $u_w\not\in S_2$ such that $(v_x,u_y)(v_x,u_w)\in E(K_n\square K_m)$ and $(v_x,u_w)\not\in N[S]$. Thus, $(v_x,u_y)$ has a private neighbour with respect to $S$, and $S$ is digitally convex. There are $(2^n-2)(2^m-2)$ such sets $S$. 

Any set of vertices containing a set of type $\{(v_1,u_{i_1}),(v_2,u_{i_2}), \dots, (v_n,u_{i_n})\}$, where each $i_j\in\{1,2,\dots,m\}$, is a dominating set in $K_n\square K_m$. Similarly, any set of vertices containing a set of type $\{(v_{j_1},u_1), (v_{j_2},u_2),\dots,(v_{j_m},u_m)\}$, where each $j_k\in\{1,2,\dots,n\}$, is a dominating set. Therefore, the only digitally convex set containing any of these sets of vertices is $V(K_n\square K_m)$. 

Two vertices $(v_x,u_y)$ and $(v_w,u_z)$ dominate the neighbourhoods of both of the vertices $(v_x,u_z)$ and $(v_w,u_y)$. So any digitally convex set containing the former pair of vertices must also contain the latter pair. Therefore, every nonempty digitally convex set in $K_n \square K_m$ must be $V(K_n\square K_m)$ or must take on the form $S_1 \times S_2$, where $S_1$ and $S_2$ are defined as above. 

Therefore, along with the empty set, the graph $K_n\square K_m$ has a total of $2+(2^n-2)(2^m-2)$ digitally convex sets. \qed

We turn now to the Cartesian product of paths, beginning with $P_n \square P_2$. As with the Cartesian product of complete graphs, many of the digitally convex sets in the product of paths are no longer digitally convex when restricted to one of the constituent paths. Thus, there is no obvious method of using the digitally convex sets of the constituent graphs to generate those of the product. We can, however, use the digitally convex sets of the graphs $P_{n-1}\square P_2$, $P_{n-2}\square P_2$ and $P_{n-3}\square P_2$ to determine those of $P_n\square P_2$. 

\begin{thm}
Let $P_n$ be the path of order $n$. Then $n_\mathscr{D}(P_1\square P_2) = 2$, $n_\mathscr{D}(P_2\square P_2) = 6$ and $n_\mathscr{D}(P_3\square P_2) = 16$ and, for $n\geq 4$, \begin{center}$n_\mathscr{D}(P_n\square P_2) = n_\mathscr{D}(P_{n-1}\square P_2)+3n_\mathscr{D}(P_{n-2}\square P_2)+2n_\mathscr{D}(P_{n-3}\square P_2).$\end{center}
\label{thm:pnp2recurrence}
\end{thm}

\noindent
\emph{Proof.} We begin by denoting the vertices of $P_n\square P_2$ by $v_1,v_2\dots,v_n$,$u_1,u_2,\dots,u_n$, with $v_iv_{i+1}\in E(P_n\square P_2)$ and $u_iu_{i+1}\in E(P_n\square P_2)$ for $i=1,2,\dots,n-1$ and $v_ju_j\in E(P_n\square P_2)$ for $j=1,2,\dots,n$. 

Now, we prove the initial conditions. Since $P_1\square P_2 \cong K_2$, we have $n_\mathscr{D}(P_1\square P_2) = n_\mathscr{D}(K_2) = 2$. Similarly, $P_2\square P_2 \cong C_4$, so $n_\mathscr{D}(P_2\square P_2) = n_\mathscr{D}(C_4) = 6$. Finally, the 16 digitally convex sets of $P_3 \square P_2$ are $\emptyset$, $\{v_1\}$, $\{v_2\}$, $\{v_3\}$, $\{u_1\}$, $\{u_2\}$, $\{u_3\}$, $\{v_1,v_3\}$, $\{v_1,u_1\}$, $\{v_3,u_3\}$, $\{u_1,u_3\}$, $\{v_1, v_3,u_1\}$, $\{v_1,u_1,u_2\}$, $\{v_2,v_3,u_3\}$, $\{v_3,u_2,u_3\}$, and $V(P_3\square P_2)$. So $n_\mathscr{D}(P_3\square P_2)=16$.

Suppose $n\geq 4$. We begin by showing $n_\mathscr{D}(P_n\square P_2) \geq n_\mathscr{D}(P_{n-1}\square P_2)+3n_\mathscr{D}(P_{n-2}\square P_2)+2n_\mathscr{D}(P_{n-3}\square P_2)$. We now construct three pairwise disjoint families $\mathscr{D}_i$, $i=1,2,3$, of digitally convex sets in $\mathscr{D}(P_n\square P_2)$ such that $|\mathscr{D}_i| = c_in_\mathscr{D}(P_{n-i}\square P_2)$, where $c_1=1,c_2=3,c_3=2$. 

To construct $\mathscr{D}_1$, let $S\in \mathscr{D}(P_{n-1} \square P_2)$. If $v_{n-1}, u_{n-1}\not\in S$, then $S$ is digitally convex in $P_n \square P_2$, because the vertices $v_n$ and $u_n$ are each a private neighbour for themselves with respect to $S$. Then, we add $S$ to $\mathscr{D}_1$. If $v_{n-1}\in S$ or $u_{n-1}\in S$, then $S\cup\{v_n,u_n\}$ is digitally convex in $P_n \square P_2$, because each vertex in $V(P_{n-1}\square P_2)-S$ must have a private neighbour with respect to $S$ in $V(P_{n-1}\square P_2)-\{v_{n-1},u_{n-1}\}$, which is also a private neighbour with respect to $S\cup\{v_n,u_n\}$ in $P_n\square P_2$. Then, we add $S\cup\{v_n,u_n\}$ to $\mathscr{D}_1$. Note that $|\mathscr{D}_1| = n_\mathscr{D}(P_{n-1}\square P_2)$, as desired. 

To construct $\mathscr{D}_2$, let $S\in \mathscr{D}(P_{n-2}\square P_2)$. If $v_{n-2}, u_{n-2}\in S$, then $S$ is digitally convex in $P_n\square P_2$, because the vertices $u_n$ and $v_n$ are private neighbours for themselves, as well as for $v_{n-1}$ and $u_{n-1}$, with respect to $S$ in $P_n\square P_2$. Then, we add $S$ to $\mathscr{D}_2$. This set is not digitally convex in $P_{n-1}\square P_2$, as the vertices $v_{n-1}$ and $u_{n-1}$ have no private neighbours with respect to $S$. The set $S\cup\{v_{n-1}\}$ is also digitally convex in $P_n\square P_2$, because the vertex $u_n$ is a private neighbour for itself, $v_n$ and $u_{n-1}$ with respect to $S$ in $P_n\square P_2$. Then, we add $S\cup\{v_1\}$ to $\mathscr{D}_2$. Similarly, $S\cup\{u_{n-1}\}$ is digitally convex in $P_n\square P_2$, so we add it to $\mathscr{D}_2$. 

If $v_{n-2}\in S$ and $u_{n-2}\not\in S$, then $S\cup\{v_n\}$ is digitally convex in $P_n\square P_2$, because the vertex $u_{n-1}$ is a private neighbour for itself, $v_{n-1}$ and $u_n$ with respect to $S\cup\{v_n\}$ in $P_n\square P_2$. Then, we add $S\cup\{v_n\}$ to $\mathscr{D}_2$. The set $S\cup\{v_{n-1}\}$ is also digitally convex in $P_n\square P_2$, because the vertex $u_n$ is a private neighbour for itself, $v_n$ and $u_{n-1}$ with respect to $S\cup\{v_{n-1}\}$. Then, we add $S\cup\{v_{n-1}\}$ to $\mathscr{D}_2$. Similarly, $S\cup\{u_{n-1}\}$ is digitally convex in $P_n\square P_2$, so we add it to $\mathscr{D}_2$. 

If $v_{n-2}\not\in S$ and $u_{n-2}\in S$ then, by the same argument as above, $S\cup\{v_{n-1}\}$, $S\cup\{u_{n-1}\}$ and $S\cup\{u_n\}$ are all digitally convex in $P_n\square P_2$. We add them all to $\mathscr{D}_2$. 

If $v_{n-2},u_{n-2}\not\in S$, then $S\cup\{v_n\}$ is digitally convex in $P_n\square P_2$ because the vertex $u_{n-1}$ is a private neighbour for itself, $u_n$ and $v_{n-1}$ with respect to $S\cup\{v_n\}$ in $P_n\square P_2$. Then, we add $S\cup\{v_n\}$ to $\mathscr{D}_2$. Similarly, $S\cup\{u_n\}$ is digitally convex in $P_n\square P_2$, so we add it to $\mathscr{D}_2$. If, in addition, $v_{n-3},u_{n-3}\not\in S$, then $S\cup\{v_n,u_n\}$ is digitally convex in $P_n\square P_2$, because both $u_{n-2},v_{n-2}\not\in N[S\cup\{v_n,u_n\}]$ in $P_n\square P_2$. So $u_{n-2}$ and $v_{n-2}$ are private neighbours for $u_{n-1}$ and $v_{n-1}$ with respect to $S\cup\{v_n,u_n\}$. Then, we add $S\cup\{v_n,u_n\}$ to $\mathscr{D}_2$. If $v_{n-3}\in S$ and $u_{n-3}\not\in S$, then $S\cup\{v_{n-1}\}$ is digitally convex in $P_n\square P_2$, because the vertex $u_{n-2}$ is a private neighbour for itself and for $v_{n-2}$, and the vertex $u_n$ is a private neighbour for itself, $v_n$ and $u_{n-1}$ with respect to $S\cup\{v_{n-1}\}$ in $P_n \square P_2$. Then, we add $S\cup\{v_{n-1}\}$ to $\mathscr{D}_2$. Similarly, if $v_{n-3}\not\in S$ and $u_{n-3}\in S$, then $S\cup\{u_{n-1}\}$ is digitally convex in $P_n\square P_2$. So we add it to $\mathscr{D}_2$. Now, we have $|\mathscr{D}_2| = 3n_\mathscr{D}(P_{n-2}\square P_2)$, as desired.  

Finally, to construct $\mathscr{D}_3$, let $S\in \mathscr{D}(P_{n-3} \square P_2)$. If $v_{n-3}, u_{n-3}\not\in S$, then $S\cup\{v_{n-1}\}$ is digitally convex in $P_n\square P_2$, because the vertex $u_{n-2}$ is a private neighbour for itself, $v_{n-2}$ and $u_{n-1}$, and the vertex $u_n$ is a private neighbour for itself and $v_n$ with respect to $S\cup\{v_{n-1}\}$. Then, we add $S\cup\{v_{n-1}\}$ to $\mathscr{D}_3$. Similarly, $S\cup\{u_{n-1}\}$ is digitally convex in $P_n\square P_2$, so we add it to $\mathscr{D}_3$. 

If $v_{n-3}\in S$ or $u_{n-3}\in S$, then $S\cup\{v_{n-2},v_n\}$ is digitally convex in $P_n\square P_2$, because the vertex $u_{n-1}$ is a private neighbour for itself, $u_{n-2}$, $v_{n-1}$ and $u_n$ with respect to $S\cup\{v_{n-2},v_n\}$ in $P_n\square P_2$. Then, we add $S\cup\{v_{n-2},v_n\}$ to $\mathscr{D}_3$. Similarly, $S\cup\{u_{n-2},u_n\}$ is digitally convex in $P_n\square P_2$, so we add it to $\mathscr{D}_3$. Now, we have $|\mathscr{D}_3| = 2n_\mathscr{D}(P_{n-3}\square P_2)$, as desired.  

Now, we have $\mathscr{D}_i \cap \mathscr{D}_j=\emptyset$ for $i\neq j$, and each $\mathscr{D}_i$, $i=1,2,3$, is a subset of $\mathscr{D}(P_n\square P_2)$. Thus $n_\mathscr{D}(P_n\square P_2) \geq |\mathscr{D}_1| + |\mathscr{D}_2| + |\mathscr{D}_3| = n_\mathscr{D}(P_{n-1}\square P_2)+3n_\mathscr{D}(P_{n-2}\square P_2)+2n_\mathscr{D}(P_{n-3}\square P_2)$.

Now, to show the reverse inequality, let $S\in \mathscr{D}(P_n\square P_2)$. \begin{enumerate}
\item[(a)] Suppose $v_n,u_n\in S$. If $v_{n-1}\in S$ or  $u_{n-1}\in S$, then each $x\not\in S$ has a private neighbour with respect to $S$ in $V(P_{n-1}\square P_2)$. Thus, $S-\{v_n,u_n\}$ is digitally convex in $P_{n-1}\square P_2$. If $v_{n-1}, u_{n-1}\not\in S$, then $v_{n-2}, u_{n-2}\not\in N[S]$. Thus, $v_{n-3}, u_{n-3}\not\in S$ and $S-\{v_n,u_n\}$ is digitally convex in $P_{n-2}\square P_2$. 

\item[(b)] Suppose $v_n\in S$ and $u_n\not\in S$. Then, $u_{n-1}\not\in N[S]$, so $v_{n-1}, u_{n-1}, u_{n-2}\not\in S$. If $v_{n-2}\in S$ and $v_{n-3}, u_{n-3}\not\in S$, then it must be the case that $u_{n-3}\not\in N[S]$. So $S-\{v_n\}$ is digitally convex in $P_{n-2}\square P_2$. If $v_{n-2}\in S$ and $v_{n-3}\in S$ or $u_{n-3}\in S$, then $S-\{v_n,v_{n-2}\}$ is digitally convex in $P_{n-3}\square P_2$. If $v_{n-2}\not\in S$, then at most one of $v_{n-3}$ and $u_{n-3}$ can be in $S$. So either $v_{n-2}\not\in N[S]$ or $u_{n-2}\not\in N[S]$. Then, $S-\{v_n\}$ is digitally convex in $P_{n-2}\square P_2$. 

\item[(c)] Suppose $v_n\not\in S$ and $u_n\in S$, then $v_{n-1}\not\in N[S]$, so $v_{n-1}, u_{n-1}, v_{n-2}\not\in S$. If $u_{n-2}\in S$ and $v_{n-3},u_{n-3}\not\in S$, then $S-\{u_n\}$ is digitally convex in $P_{n-2}\square P_2$. If $u_{n-2}\in S$ and $v_{n-3}\in S$ or $u_{n-3}\in S$, then $S-\{u_n,u_{n-2}\}$ is digitally convex in $P_{n-3}\square P_2$. If $u_{n-2}\not\in S$, then $S-\{u_n\}$ is digitally convex in $P_{n-2}\square P_2$. 

\item[(d)] Suppose $v_n,u_n\not\in S$. Then, at most one of $v_{n-1}$ and $u_{n-1}$ can be in $S$. If $v_{n-1}\in S$ and at least one of $v_{n-2}$ and $u_{n-2}$ is in $S$, then $S-\{v_{n-1}\}$ is digitally convex in $P_{n-2}\square P_2$. If $v_{n-1}\in S$ and $v_{n-2}, u_{n-2}\not\in S$, then it must be the case that $u_{n-2}\not\in N[S]$. So $u_{n-3}\not\in S$. If $v_{n-3}\in S$, then $S-\{v_{n-1}\}$ is digitally convex in $P_{n-2}\square P_2$. If $v_{n-3}\not\in S$, then $S-\{v_{n-1}\}$ is digitally convex in $P_{n-3}\square P_2$. Similarly, if $u_{n-1}\in S$ and at least one of $v_{n-2}$ and $u_{n-2}$ is in $S$, then $S-\{u_{n-1}\}$ is digitally convex in $P_{n-2}\square P_2$. If $u_{n-1}\in S$, $v_{n-2}, u_{n-2}\not\in S$ and $u_{n-3}\in S$, then $S-\{u_{n-1}\}$ is digitally convex in $P_{n-2}\square P_2$. If $u_{n-1}\in S$, $v_{n-2}, u_{n-2}\not\in S$ and $u_{n-3}\not\in S$, then $S-\{u_{n-1}\}$ is digitally convex in $P_{n-3}\square P_2$.

\end{enumerate}

Each digitally convex set in $P_{n-1}\square P_2$ has been counted here at most once, each digitally convex set in $P_{n-2}\square P_2$ at most three times, and each digitally convex set in $P_{n-3}\square P_2$ at most twice.
Refer to Table~\ref{table:1}
for a summary of which digitally convex sets in $P_{n-2}\square P_2$ and $P_{n-3}\square P_2$ are counted in each part of the above argument. 
Therefore, $n_\mathscr{D}(P_n\square P_2) \leq n_\mathscr{D}(P_{n-1}\square P_2)+3n_\mathscr{D}(P_{n-2}\square P_2)+2n_\mathscr{D}(P_{n-3}\square P_2)$. \qed

\begin{table}[ht]
\centering
\begin{tabular}{ |c|c| } 
\hline
\multicolumn{2}{|c|}{$P_{n-2}\square P_2$} \\
 \hline
 \begin{tikzpicture}[thick, every node/.style={circle, draw=black, fill=black, inner sep=2}]
 
 \node[label=left:$v_{n-2}$] (n0) at (0,0){};
 \node[label=right:$u_{n-2}$] (n1) at (1,0){};
 \node[label=left:$v_{n-3}$] (n2) at (0,1){};
 \node[label=right:$u_{n-3}$] (n3) at (1,1){};
 \node[draw=none, fill=none] (n4) at (0.5,1.6){};
 
 \draw (n0) -- (n1);
 \draw (n0) -- (n2);
 \draw (n1) -- (n3);
 \draw (n2) -- (n3);
 \draw (n2) -- (0, 1.35);
 \draw (n3) -- (1, 1.35);
 \end{tikzpicture} & (a), (b), and (c)  \\ 
 \hline
 \begin{tikzpicture}[thick, every node/.style={circle, draw=black, fill=black, inner sep=2}]
 
 \node[fill=none, label=left:$v_{n-2}$] (n0) at (0,0){};
 \node[label=right:$u_{n-2}$] (n1) at (1,0){};
 \node[label=left:$v_{n-3}$] (n2) at (0,1){};
 \node[label=right:$u_{n-3}$] (n3) at (1,1){};
 \node[draw=none, fill=none] (n4) at (0.5,1.6){};
 
 \draw (n0) -- (n1);
 \draw (n0) -- (n2);
 \draw (n1) -- (n3);
 \draw (n2) -- (n3);
 \draw (n2) -- (0, 1.35);
 \draw (n3) -- (1, 1.35);
 \end{tikzpicture}& (b), and twice in (d)  \\ 
 \hline
 \begin{tikzpicture}[thick, every node/.style={circle, draw=black, fill=black, inner sep=2}]
 
 \node[label=left:$v_{n-2}$] (n0) at (0,0){};
 \node[fill=none, label=right:$u_{n-2}$] (n1) at (1,0){};
 \node[label=left:$v_{n-3}$] (n2) at (0,1){};
 \node[label=right:$u_{n-3}$] (n3) at (1,1){};
 \node[draw=none, fill=none] (n4) at (0.5,1.6){};
 
 \draw (n0) -- (n1);
 \draw (n0) -- (n2);
 \draw (n1) -- (n3);
 \draw (n2) -- (n3);
 \draw (n2) -- (0, 1.35);
 \draw (n3) -- (1, 1.35);
 \end{tikzpicture}& (c), and twice in (d)  \\ 
 \hline
\begin{tikzpicture}[thick, every node/.style={circle, draw=black, fill=black, inner sep=2}]
 
 \node[label=left:$v_{n-2}$] (n0) at (0,-0.75){};
 \node[label=right:$u_{n-2}$] (n1) at (1,-0.75){};
 \node[fill=none, label=left:$v_{n-3}$] (n2) at (0,0.25){};
 \node[label=right:$u_{n-3}$] (n3) at (1,0.25){};
 \node[draw=none, fill=none] (n4) at (0.5,1){};
 
 \draw (n0) -- (n1);
 \draw (n0) -- (n2);
 \draw (n1) -- (n3);
 \draw (n2) -- (n3);
 \draw (n2) -- (0, 0.6);
 \draw (n3) -- (1, 0.6);
 \end{tikzpicture} \begin{tikzpicture}[thick, every node/.style={circle, draw=none, fill=none, inner sep=2}]
 \node (n0) at (0,1){and};
 \end{tikzpicture} 
 \begin{tikzpicture}[thick, every node/.style={circle, draw=black, fill=black, inner sep=2}]
 
 \node[label=left:$v_{n-2}$] (n0) at (0,-0.75){};
 \node[label=right:$u_{n-2}$] (n1) at (1,-0.75){};
 \node[label=left:$v_{n-3}$] (n2) at (0,0.25){};
 \node[fill=none, label=right:$u_{n-3}$] (n3) at (1,0.25){};
 \node[draw=none, fill=none] (n4) at (0.5,1){};
 
 \draw (n0) -- (n1);
 \draw (n0) -- (n2);
 \draw (n1) -- (n3);
 \draw (n2) -- (n3);
 \draw (n2) -- (0, 0.6);
 \draw (n3) -- (1, 0.6);
 \end{tikzpicture} & (b), (c), and (d)  \\ 
 \hline
 \multicolumn{2}{|c|}{$P_{n-3}\square P_2$} \\
 \hline
\begin{tikzpicture}[thick, every node/.style={circle, draw=black, fill=black, inner sep=2}]

 \node[fill=none, label=below:$v_{n-3}$] (n2) at (0,1){};
 \node[label=below:$u_{n-3}$] (n3) at (1,1){};
 \node[draw=none, fill=none] (n4) at (0.5,1.6){};

 \draw (n2) -- (n3);
 \draw (n2) -- (0, 1.35);
 \draw (n3) -- (1, 1.35);
 \end{tikzpicture} \begin{tikzpicture}[thick, every node/.style={circle, draw=none, fill=none, inner sep=2}]
 \node (n0) at (0,1){and};
 \end{tikzpicture} 
 \begin{tikzpicture}[thick, every node/.style={circle, draw=black, fill=black, inner sep=2}]
 
 \node[label=below:$v_{n-3}$] (n2) at (0,1){};
 \node[fill=none, label=below:$u_{n-3}$] (n3) at (1,1){};
 \node[draw=none, fill=none] (n4) at (0.5,1.6){};

 \draw (n2) -- (n3);
 \draw (n2) -- (0, 1.35);
 \draw (n3) -- (1, 1.35);
 \end{tikzpicture} \begin{tikzpicture}[thick, every node/.style={circle, draw=none, fill=none, inner sep=2}]
 \node (n0) at (0,1){and};
 \end{tikzpicture}
 \begin{tikzpicture}[thick, every node/.style={circle, draw=black, fill=black, inner sep=2}]
 
 \node[fill=none, label=below:$v_{n-3}$] (n2) at (0,1){};
 \node[fill=none, label=below:$u_{n-3}$] (n3) at (1,1){};
 \node[draw=none, fill=none] (n4) at (0.5,1.6){};

 \draw (n2) -- (n3);
 \draw (n2) -- (0, 1.35);
 \draw (n3) -- (1, 1.35);
 \end{tikzpicture} & (b) and (c)  \\ 
 \hline
\begin{tikzpicture}[thick, every node/.style={circle, draw=black, fill=black, inner sep=2}]
 
 \node[label=left:$v_{n-3}$] (n2) at (0,1){};
 \node[label=right:$u_{n-3}$] (n3) at (1,1){};
 \node[draw=none, fill=none] (n4) at (0.5,1.6){};

 \draw (n2) -- (n3);
 \draw (n2) -- (0, 1.35);
 \draw (n3) -- (1, 1.35);
 \end{tikzpicture} & twice in (d)  \\ 
 \hline
\end{tabular}
\caption{A summary of the counting argument in Theorem~\ref{thm:pnp2recurrence}, with white vertices in $S$}
\label{table:1}
\end{table}

Note that, in addition to proving the given recurrence, this proof of Theorem~\ref{thm:pnp2recurrence} provides a method of generating the collection of digitally convex sets of $P_n \square P_2$ from those of $P_{n-1} \square P_2$, $P_{n-2}\square P_2$ and $P_{n-3}\square P_2$. 

Given a set $\mathscr{A}$ of $n\times m$ binary arrays, we let $\mathscr{A}^*$ be the set of arrays obtained as follows. For each array $A$ in $\mathscr{A}$, construct a new array $A^*$ by taking the minimum value of corresponding elements of $A$ and their horizontal and vertical neighbours. In other words, each element of $A^*$ is the minimum value over the closed neighbourhood of the corresponding element of $A$. 

As an example, consider the following array $A$
$$\begin{array}{ccc}
1 & 1 & 0 \\
1 & 1 & 1 \\
0 & 1 & 1
\end{array}
$$

Taking the minimum value over the closed neighbourhood of each element in the array produces the array $A^*$:
$$\begin{array}{ccc}
1 & 0 & 0\\
0 & 1 & 0 \\
0 & 0 & 1
\end{array}$$

Note that in this process, two distinct arrays $A_1$ and $A_2$ can produce the same array $A^*$. We now establish a bijection between these $n\times m$ binary arrays and the digitally convex sets in $\mathscr{D}(P_n \square P_m)$. 

\begin{thm}
Let $\mathscr{A}_{n,m}$ be the set of all $n\times m$ binary arrays. Then, $n_\mathscr{D}(P_n\square P_m)$ = $|\mathscr{A}_{n,m}^*|$. 
\label{thm:pnpm}
\end{thm}

\noindent
\emph{Proof.} First, we label the vertices of the product $P_n \square P_m$. Let the vertices of $P_n$ be $u_1,u_2,\dots,u_n$, with $u_i u_{i+1}\in E(P_n)$ for $i=1,2,\dots,n-1$, and let the vertices of $P_m$ be $v_1,v_2,\dots,v_m$, with $v_jv_{j+1}\in E(P_m)$ for $j=1,2,\dots,m-1$. Then, the vertices of $P_n\square P_m$ have the form $(u_i,v_j)$. 

Now we show a bijection between the digitally convex sets in $\mathscr{D}(P_n\square P_m)$ and the arrays in $\mathscr{A}_{n,m}^*$. Let $A^* \in \mathscr{A}_{n,m}^*$ and consider the set $S = \{(u_i,v_j)\mid a^*_{i,j}=1\}$. Each vertex $(u_x,v_y)\not\in S$ corresponds to an entry $a^*_{x,y}$ that has value 0 in $A^*$. Then, either the corresponding entry in $A$ also has value 0, or it has value 1 and has a horizontal or vertical neighbour with value 0. In the first case, every entry in the closed neighbourhood of $a^*_{x,y}$ also has value 0 in $A^*$. In $P_n \square P_m$, this means that none of the vertices in $N[(u_x,v_y)]$ is in $S$, so $(u_x,v_y)$ is its own private neighbour. In the second case, there is an entry $a_{w,z}$ in the closed neighbourhood of $a_{x,y}$ which has value 0 in $A$. Then, in $A^*$, every entry in the closed neighbourhood of $a^*_{w,z}$ has value 0, including $a^*_{x,y}$. In $P_n\square P_m$, this means that none of the vertices in $N[(u_w,v_z)]$ is in $S$ and $(u_w,v_z)(u_x,v_y)\in E(P_n\square P_m)$, so the vertex $(u_w,v_z)$ is a private neighbour for $(u_x,v_y)$ with respect to $S$. Therefore, $S$ is digitally convex in $P_n \square P_m$. 

It is clear from the construction of $S$ that this mapping from $\mathscr{A}^*_{n,m}$ to $\mathscr{D}(P_n\square P_m)$ is injective. It remains to be shown that the mapping is surjective. Consider $S\in \mathscr{D}(P_n\square P_m)$ and let $B$ be the $n\times m$ array with $b_{i,j}=1$ if $(u_i,v_j)\in S$ and $b_{i,j}=0$ otherwise. Then, let $C$ be the $n\times m$ array whose entries are the maximum over the closed neighbourhood of the corresponding entry in $B$. In other words, $c_{i,j}=1$ if any of the entries in the closed neighbourhood of $b_{i,j}$ has value 1, and $c_{i,j}=0$ otherwise. Clearly, $C\in\mathscr{A}$ and now we show that $C^*=B$. By construction of $C$, each entry of $C^*$ whose corresponding entry in $B$ has value 1 also has value 1 in $C^*$. So if $C^* \neq B$, then there is some $i,j$ with $c^*_{i,j}=1$ and $b_{i,j}=0$. This means that, in $C$, each entry in the closed neighbourhood of $c_{i,j}$ has value 1. However, the entries in $C$ are defined to be 1 because their corresponding entry in $B$ has a 1 in its closed neighbourhood. In other words, the entries in the closed neighbourhood of $b_{i,j}$ each either have value 1 or have a horizontal or vertical neighbour with value 1. In terms of the set $S$, this corresponds to a vertex $(u_i,v_j)$ with every vertex in $N[(u_i,v_j)]$ in $N[S]$, i.e.~$(u_i,v_j)$ has no private neighbour with respect to $S$ in $P_n\square P_m$. This contradicts $S$ being digitally convex and thus $C^*=B$. 

It is clear that $B$ gets mapped to the digitally convex set $S$, using the mapping described above. Therefore, $n_\mathscr{D}(P_n\square P_m) = |\mathscr{A}^*_{n,m}|$. \qed

The number of digitally convex sets of $P_n\square P_m$ follows the OEIS sequence A217637~\cite{oeis1}. The OEIS notes an observation from Andrew Howroyd that this sequence also enumerates the maximal independent sets in the graph $P_n \square P_m \square P_2$. Euler, Oleksik and Skupie\'n~\cite{euler} prove this equivalence for $m=2$ and for $m=3$. However, the correspondence between the maximal independent sets in $P_n\square P_m\square P_2$ and the digitally convex sets in $P_n\square P_m$ is not clear, even for very small values of $n$ and $m$. 

\section{Conclusion}

In this paper, we established a linear recurrence that is satisfied by the number of cyclic binary $n$-bit strings whose blocks each have length at least $k$, for some $k\geq 2$. We then showed that these cyclic binary strings can be used to enumerate the digitally convex sets of powers of cycles. It is possible that cyclic binary strings with other properties can be used to enumerate the digitally convex sets of graphs that can be constructed by adding additional edges to a cycle, $C_n$. 

We provided a closed formula for the number of digitally convex sets in the graph $K_n\square K_m$. We established a recurrence relation for $n_\mathscr{D}(P_n\square P_2)$ and defined a class of $n\times m$ binary arrays for which there is a one-to-one correspondence with the digitally convex sets of $P_n\square P_m$. The problem of enumerating the digitally convex sets for other graph products, such as the strong product, the direct product and the categorical product remains open.

\bibliography{refs}
\end{document}